\documentclass[12pt]{amsart}
\usepackage{amssymb,latexsym}
\usepackage{pdfsync}

\newdimen\AAdi%
\newbox\AAbo%
\def\AAk#1#2{\s_etbox\AAbo=\hbox{#2}\AAdi=\wd\AAbo\kern#1\AAdi{}}%
\def\AAr#1#2#3{\s_etbox\AAbo=\hbox{#2}\AAdi=\ht\AAbo\raise#1\AAdi\hbox{#3}}%
\font\tenmsb=msbm10 at 12pt
\font\sevenmsb=msbm7 at 8pt
\font\fivemsb=msbm5 at 6pt
\newfam\msbfam
\textfont\msbfam=\tenmsb
\scriptfont\msbfam=\sevenmsb
\scriptscriptfont\msbfam=\fivemsb
\def\Bbb#1{{\tenmsb\fam\msbfam#1}}

\newcommand{\beq}{\begin{equation}}
\newcommand{\eeq}{\end{equation}}
\newcommand{\beqr}{\begin{eqnarray}}
\newcommand{\eeqr}{\end{eqnarray}}
\newcommand{\ba}{\begin{array}}
\newcommand{\ea}{\end{array}}

\makeatletter

\newcommand{\Rmnum}[1]{\expandafter\@slowromancap\romannumeral #1@}
\makeatother

\begin{document}

\newtheorem{thm}{Theorem}
\newtheorem{lem}{Lemma}
\newtheorem{cor}{Corollary}
\newtheorem{rem}{Remark}
\newtheorem{pro}{Proposition}
\newtheorem{defi}{Definition}
\newtheorem{eg}{Example}
\newtheorem*{claim}{Claim}
\newtheorem{conj}[thm]{Conjecture}
\newcommand{\noi}{\noindent}
\newcommand{\dis}{\displaystyle}
\newcommand{\mint}{-\!\!\!\!\!\!\int}
\numberwithin{equation}{section}

\def \bx{\hspace{2.5mm}\rule{2.5mm}{2.5mm}}
\def \vs{\vspace*{0.2cm}}
\def\hs{\hspace*{0.6cm}}
\def \ds{\displaystyle}
\def \p{\partial}
\def \O{\Omega}
\def \o{\omega}
\def \b{\beta}
\def \m{\mu}
\def \l{\lambda}
\def\L{\Lambda}
\def \ul{u_\lambda}
\def \D{\Delta}
\def \d{\delta}
\def \k{\kappa}
\def \s{\sigma}
\def \e{\varepsilon}
\def \a{\alpha}
\def \tf{\tilde{f}}
\def\cqfd{%
\mbox{ }%
\nolinebreak%
\hfill%
\rule{2mm} {2mm}%
\medbreak%
\par%
}
\def \pr {\noindent {\it Proof.} }
\def \rmk {\noindent {\it Remark} }
\def \esp {\hspace{4mm}}
\def \dsp {\hspace{2mm}}
\def \ssp {\hspace{1mm}}

\def\la{\langle}\def\ra{\rangle}

\def \u{u_+^{p^*}}
\def \ui{(u_+)^{p^*+1}}
\def \ul{(u^k)_+^{p^*}}
\def \energy{\int_{\R^n}\u }
\def \sk{\s_k}
\def \mo{\mu_k}
\def\cal{\mathcal}
\def \I{{\cal I}}
\def \J{{\cal J}}
\def \K{{\cal K}}
\def \OM{\overline{M}}

\def\n{\nabla}

\def\fk{{{\cal F}}_k}
\def\M1{{{\cal M}}_1}
\def\Fk{{\cal F}_k}
\def\Fl{{\cal F}_l}
\def\FF{\cal F}
\def\Gk{{\Gamma_k^+}}
\def\n{\nabla}
\def\uuu{{\n ^2 u+du\otimes du-\frac {|\n u|^2} 2 g_0+S_{g_0}}}
\def\uuug{{\n ^2 u+du\otimes du-\frac {|\n u|^2} 2 g+S_{g}}}
\def\sku{\sk\left(\uuu\right)}
\def\qed{\cqfd}
\def\vvv{{\frac{\n ^2 v} v -\frac {|\n v|^2} {2v^2} g_0+S_{g_0}}}
\def\vvs{{\frac{\n ^2 \tilde v} {\tilde v}
 -\frac {|\n \tilde v|^2} {2\tilde v^2} g_{S^n}+S_{g_{S^n}}}}
\def\skv{\sk\left(\vvv\right)}
\def\tr{\hbox{tr}}
\def\pO{\partial \Omega}
\def\dist{\hbox{dist}}
\def\RR{\Bbb R}\def\R{\Bbb R}
\def\C{\Bbb C}
\def\B{\Bbb B}
\def\N{\Bbb N}
\def\Q{\Bbb Q}
\def\Z{\Bbb Z}
\def\PP{\Bbb P}
\def\EE{\Bbb E}
\def\F{\Bbb F}
\def\G{\Bbb G}
\def\H{\Bbb H}
\def\SS{\Bbb S}\def\S{\Bbb S}

\def\div{\hbox{div}\,}

\def\lcf{{locally conformally flat} }

\def\circledwedge{\setbox0=\hbox{$\bigcirc$}\relax \mathbin {\hbox
to0pt{\raise.5pt\hbox to\wd0{\hfil $\wedge$\hfil}\hss}\box0 }}

\def\sss{\frac{\s_2}{\s_1}}

\date{May 13, 2024}
\title[ Liouville theorems for ancient solutions ]{ Liouville theorems for ancient solutions to the V-harmonic map heat flows $\rm \Rmnum{2}$}

\author{}

\author[Chen]{Qun Chen}
\address{School of Mathematics and Statistics\\ Wuhan University\\Wuhan 430072,
China
 }
 \email{qunchen@whu.edu.cn}

\author[Qiu]{Hongbing Qiu$^*$}
\address{School of Mathematics and Statistics\\ Wuhan University\\Wuhan 430072,
China
 }
 \email{hbqiu@whu.edu.cn}

 \renewcommand{\thefootnote}{\fnsymbol{footnote}}

\footnotetext{\hspace*{-5mm} \begin{tabular}{@{}r@{}p{13.4cm}@{}}

 $^*$ & Corresponding author \\

\end{tabular}}

\renewcommand{\thefootnote}{\arabic{footnote}}

\begin{abstract}

When the domain is a complete noncompact Riemannian manifold with nonnegative Bakry--Emery Ricci curvature and the target is a complete Riemannian manifold with sectional curvature bounded above by a positive constant, by carrying out refined gradient estimates, we obtain a better Liouville theorem for ancient solutions to the $V$-harmonic map heat flows.  Furthermore,  we can also derive a Liouville theorem for quasi-harmonic maps under an exponential growth condition.  

\vskip12pt



\end{abstract}
\maketitle
\section{Introduction}

We call that a map $u$ from a Riemannian manifold $(M^m, g)$ to another Riemannian manifold manifold $(N^n,  h)$ is a $V$-harmonic map if it satisfies
\[
\tau_V(u):= \tau(u)+ du(V) = 0,
\]
where $V$ is a vector field on $M$, and $\tau(u)$, the tension field of the map $u$, is given by 
\[
\tau(u) = g^{\a\b} \left( \frac{\p^2 u^k}{\p x^\a \p x^\b} - \Gamma^\sigma_{\a\b}\frac{\p u^k}{\p x^\sigma}  +  \Gamma^k_{ji}\frac{\p u^i}{\p x^\a} \frac{\p u^j}{\p x^\b} \right)\frac{\p}{\p y^k}
\]
in local coordinates (we let $\{x^\alpha\}$ and $\{y^j\}$ be local coordinates on $M$ and $N$, respectively. $\Gamma^\s_{\a\b}$ and $\Gamma^k_{ji}$ are the Christoffel's symbols of $M$ and $N$, respectively). This is a generalization of the
usual harmonic map. The corresponding $V$-harmonic map heat flow is 
\begin{equation}\label{eqn-anc111}
\frac{\p}{\p t} u = \tau_V(u).
\end{equation}

The notion of $V$-harmonic maps was introduced by the first author, Jost and Wang in \cite{CheJosWan15}. It includes the Hermitian harmonic maps introduced and studied in \cite{JY93}, the Weyl harmonic maps from a Weyl manifold into a Riemannian manifold \cite{Kok09}, the affine harmonic maps mapping from an affine manifold into a Riemannian manifold \cite{JS09, JS11}, and harmonic maps from a Finsler manifold into a Riemannian manifold \cite{Cen00,HS08, Mo01,MY05, MW09}. In \cite{CheJosWan15}, the authors showed the existence of $V$-harmonic maps from compact manifolds with boundary. Later, Chen-Jost-Qiu \cite{CheJosQiu12} generalized this result to the case that the domain is a complete manifold and then the second author \cite{Qiu17} obtained the global existence of $V$-harmonic map heat flows from complete manifolds.  

When the domain is a complete noncompact Riemannian manifold with nonnegative Ricci curvature and the target is a simply connected Riemannian manifold with nonpositive sectional curvature, Wang \cite{Wan11} showed that any ancient solution to the harmonic map heat flows must be constant under certain growth condition near infinity. See also the related work \cite{CLQ22, CM21, Hua18,  KLLS23, KS21, LZ19, SZ}, etc. Recently, the authors \cite{CQ23} generalized this result of \cite{Wan11} to the case for $V$-harmonic maps. Furthermore, they also proved that when the Bakry--Emery Ricci curvature of the domain is nonnegative and the target is a complete Riemannian manifold with sectional curvature bounded above by a positive constant, then any ancient solution to (\ref{eqn-anc111}) with its image contained in a regular ball has to be a constant map \cite{CQ23}.

In this note, we continue to study the ancient solution to the $V$-harmonic map heat flow. With the aid of an appropriate auxiliary function, we obtain the following gradient estimate by using the technique of \cite{SZ}.

\begin{thm} \label{thm1}

Let $(M^m, g)$ be a complete noncompact Riemannian manifold with Bakry--Emery Ricci curvature
\[
{\rm Ric}_V := {\rm Ric} - \frac{1}{2}L_V g  \geq 0,
\]
where ${\rm Ric}$ is the Ricci curvature of $M$ and $L_V$ is the Lie derivative. Let $(N^n, h)$ be a complete Riemannian manifold with sectional curvature bounded above by a positive constant $\k$.
 Let $r, \rho$ be the distance function on $M, N$ from $\tilde p_0\in M, p_0\in N$ respectively. Let $u(x, t)$ be a solution to (\ref{eqn-anc111}). Denote $B_R(\tilde p_0):= \{ x\in M| r(x) \leq R \}, Q_{R, T}:= B_R(\tilde p_0) \times [t_0-T, t_0]\subset M\times (-\infty, +\infty)$.  
 Set 
 \[
 \varphi=1-\cos \sqrt{\k}\rho,  \quad \quad b=\frac{1}{2}\left( 1+\sup_{Q_{R, T}}\varphi\circ u \right).
\] 
  Assume that $u(Q_{R, T})$ is contained in an open geodesic ball of $N$ centered at $p_0$ of radius $\frac{\pi}{2\sqrt\k}$. Then we have
\begin{equation}\label{eqn-Pa130}
\aligned
 \sup_{Q_{\frac{R}{2}, \frac{T}{2}}} \frac{e(u)}{(b - \varphi\circ u)^{2}} 
\leq & C \left( \frac{1}{R^2}\cdot \left( \sup_{Q_{R, T}}\left( \frac{\pi}{2\sqrt \k} - \rho\circ u \right)^{-2} 
 +  \sup_{Q_{R, T}}\left( \frac{\pi}{2\sqrt \k} - \rho\circ u \right)^{-4} \right) \right. \\
& \left.  +\frac{1}{R}\cdot \sup_{Q_{R, T}}\left( \frac{\pi}{2\sqrt \k} - \rho\circ u \right)^{-2} +\frac{1}{T}\cdot \sup_{Q_{R, T}}\left( \frac{\pi}{2\sqrt \k} - \rho\circ u \right)^{-2} \right),
\endaligned
\end{equation} 
 where $C$ is a positive constant independent of $R$ and $T$.

\end{thm}

The above estimate (\ref{eqn-Pa130}) leads to a Liouville theorem for ancient solutions to the $V$-harmonic map heat flows. 

\begin{thm} \label{thm2}

Let $(M^m, g)$ be a complete noncompact Riemannian manifold with nonnegative Bakry--Emery Ricci curvature, and
$(N^n, h)$  a complete Riemannian manifold with sectional curvature bounded above by a positive constant $\k$. Let $u(x, t)$ be an ancient solution to (\ref{eqn-anc111}). If the image of $u$ is contained in an open geodesic ball of $N$ centered at $p_0$ of radius $\frac{\pi}{2\sqrt \k}$ and 
\[
\left(\frac{\pi}{2\sqrt \k} - \rho\circ u(x, t)\right)^{-1} = o\left(\sqrt{ r(x)+\sqrt{ |t| } } \right)
\]
 near infinity, 
then $u$ is a constant.

\end{thm}

\begin{rem}

(1) The authors \cite{CQ23} showed that that such a Liouville theorem for ancient solutions $u(x, t)$ to the $V$-harmonic map heat flows holds under the condition that the image of $u$ is contained in a regular ball of $N^n$. Namely, $\rho\circ u$ is bounded by some constant $\widetilde R< \frac{\pi}{2\sqrt{\k}}$.

(2) Kunikawa--Sakurai \cite{KS21} obtained a similar result of a solution to the backward harmonic heat flow along a backward super Ricci flow.

\end{rem}

Consequently,  we derive a Liouville theorem for $V$-harmonic maps as follows.

\begin{cor}\label{thm-main4}

Let $(M^m,  g)$ be an $m$-dimensional complete noncompact Riemannian manifold with nonnegative Bakry-Emery Ricci curvature, $(N^n,  h)$ a complete Riemannian manifold with sectional curvature bounded above by a positive constant $\k$. Let $u: M\to N$ be a $V$-harmonic map. If $u(M)$ is contained in an open geodesic ball of $N$ centered at $p_0$ of radius $\frac{\pi}{2\sqrt{\k}}$, and 
\[
\left( \frac{\pi}{2\sqrt{\k}} - \rho\circ u \right)^{-1} = o(R^{\frac{1}{2}}).
\]
where $\rho$ denotes the distance on $N$ from $p_0$. Then $u$ must be a constant map.
\end{cor}

\begin{rem}

The second author \cite{Qiu17} proved that such a Liouville theorem for $V$-harmonic maps holds under the condition that the image of $M^m$ under the map $u$ is contained in a regular ball of $N^n$. Namely, $\rho\circ u$ is bounded by some constant $R_0< \frac{\pi}{2\sqrt{\k}}$. Hence the condition in Corollary \ref{thm-main4} is weaker than the one in \cite{Qiu17}. 

\end{rem}

In particular, when $u: (\mathbb R^m, g_0) \to (N^n,  h)$ is a quasi-harmonic map, that is, $u$ satisfies
\begin{equation}\label{eqn-QH1}\aligned
\tau(u) - \frac{1}{2}\la x, \n u \ra=0,
\endaligned
\end{equation}
where $g_0:= \sum_{i=1}^m dx_i^{2}$ is the metric of the Euclidean space $\mathbb R^m$ and $x=(x_1, \cdot\cdot\cdot, x_m)$. 
We could obtain the following Liouville theorem for quasi-harmonic maps.

\begin{thm}\label{thm-QH}

Let $(N^n,  h)$ be a complete Riemannian manifold with sectional curvature bounded above by a positive constant $\k$. Let $u: (\mathbb R^m, g_0) \to (N^n,  h)$ be a quasi-harmonic map. Suppose that $u(M)$ is contained an open geodesic ball of $N$ centered at $p_0$ of radius $\frac{\pi}{2\sqrt \k}$. Let $\e_0\in (0, \frac{1}{4})$ be a constant which is sufficiently small. If for any constant $\a>0$, 
\begin{equation*}\label{eqn-thm-QH}
\left( \frac{\pi}{2\sqrt \k} - \rho\circ u \right)^{-1} \leq Ce^{\a (\frac{1}{4}-\e_0)|x|^2},
\end{equation*}
where $C$ is a positive constant and $\rho$ denotes the distance on $N$ from $p_0$. Then $u$ must be a constant map.

\end{thm}

\begin{rem}

When the target manifold $N$ is a complete Riemannian manifold with sectional curvature bounded above by a positive constant $\kappa$, Li-Wang \cite{LW} showed that the quasi-harmonic map from $\mathbb R^m$ into a regular ball of $N$ must be a constant map. That is, $\rho\circ u$ is bounded by some constant $\rho_0 < \frac{\pi}{2\sqrt{\kappa}}$. Thus Theorem \ref{thm-QH} generalizes their result.

\end{rem}

\vskip24pt

\section{Proofs of main theorems}

Let $\D_V := \D + \la V, \n \cdot \ra$, where $\Delta$ is the Laplace--Beltrami operator on $M$.

\vskip8pt

\noindent{\bf Proof of Theorem \ref{thm1}}
Multiplying the metric tensor by a suitable constant we may assume that
the upper bound of the sectional curvature of $N$ is 1.
The Hessian comparison theorem implies 
\begin{equation}\label{eqn-Para00}
{\rm Hess}(\rho) \geq  \cot\rho\cdot (h - d\rho \otimes d\rho).
\end{equation}
It follows that
\begin{equation}\label{eqn-Para1phi}\aligned
{\rm Hess}(\varphi) = & \varphi' {\rm Hess} (\rho) + \varphi'' d\rho\otimes d\rho \\
=& \sin\rho{\rm Hess} (\rho) + \cos\rho d\rho\otimes d\rho \\
\geq & \sin\rho \cot\rho\cdot(h - d\rho \otimes d\rho) + \cos\rho d\rho\otimes d\rho \\
= &( \cos\rho) h.
\endaligned
\end{equation}

Let
\[
f:= \frac{e(u)}{(b-\varphi\circ u)^2}.
\]
Direct computation gives us
\begin{equation}\label{eqn-Pa1}\aligned
\left( \D_V -\p_t \right) f = & \frac{\left( \D_V -\p_t \right) e(u)}{(b-\varphi\circ u)^2} + \frac{2e(u)\left(\D_V -\p_t\right)(\varphi\circ u)}{(b-\varphi\circ u)^3} \\
& + \frac{4\la \n e(u), \n(\varphi\circ u) \ra}{(b-\varphi\circ u)^3} + \frac{6e(u)|\n(\varphi\circ u)|^2}{(b-\varphi\circ u)^4}.
\endaligned
\end{equation}
Since
\begin{equation}\label{eqn-Pa2}\aligned
\frac{\la \n f, \n(\varphi\circ u) \ra}{b-\varphi\circ u} =  \frac{\la \n e(u), \n(\varphi\circ u) \ra}{(b-\varphi\circ u)^3} + \frac{2e(u)|\n(\varphi\circ u)|^2}{(b-\varphi\circ u)^4}.
\endaligned
\end{equation}
From (\ref{eqn-Pa1}) and (\ref{eqn-Pa2}), we derive
\begin{equation}\label{eqn-Pa3}\aligned
\left( \D_V -\p_t \right) f = & \frac{\left( \D_V -\p_t \right) e(u)}{(b-\varphi\circ u)^2} + \frac{2e(u)\left(\D_V -\p_t\right)(\varphi\circ u)}{(b-\varphi\circ u)^3} \\
& + \frac{2\la \n f, \n(\varphi\circ u) \ra}{b-\varphi\circ u} + \frac{2e(u)|\n(\varphi\circ u)|^2}{(b-\varphi\circ u)^4} \\
& +  \frac{2\la \n e(u), \n(\varphi\circ u) \ra}{(b-\varphi\circ u)^3}.
\endaligned
\end{equation}
The composition formula,  (\ref{eqn-anc111}) and (\ref{eqn-Para1phi}) imply 
\begin{equation}\label{eqn-Pa4}\aligned
\left(\D_V - \p_t \right){\varphi\circ u} = &\sum_{i=1}^m {\rm Hess} (\varphi) (du(e_i), du(e_i))\circ u 
+ dh(\tau_V(u)-\p_t u) \\
=& \sum_{i=1}^m {\rm Hess} (\varphi) (du(e_i), du(e_i)) \circ u \\
\geq & \cos(\rho\circ u)e(u).  
\endaligned
\end{equation}
By Lemma 1 in \cite{Qiu17}, we get
\begin{equation}\label{eqn-Pa5}\aligned
\left(\D_V - \p_t \right) e(u) \geq 2|\n du|^2 - 2 e(u)^2.
\endaligned
\end{equation}
Substituting (\ref{eqn-Pa4}) and (\ref{eqn-Pa5}) into (\ref{eqn-Pa3}), we have
\begin{equation}\label{eqn-Pa6}\aligned
\left(\D_V - \p_t \right) f \geq & 2(b - \varphi\circ u)^2 \left( \frac{\cos(\rho\circ u)}{b-\varphi\circ u}-1 \right) f^2 + \frac{2|\n du|^2}{(b-\varphi\circ u)^2} \\
&+ \frac{2\la \n f, \n(\varphi\circ u) \ra}{b-\varphi\circ u} 
+ \frac{2e(u)|\n(\varphi\circ u)|^2}{(b-\varphi\circ u)^4} +  \frac{2\la \n e(u), \n(\varphi\circ u) \ra}{(b-\varphi\circ u)^3}.
\endaligned
\end{equation}
By the Schwartz inequality, we obtain
\begin{equation}\label{eqn-Pa7}\aligned
\left| \frac{2\la \n e(u), \n(\varphi\circ u) \ra}{(b-\varphi\circ u)^3} \right| \leq \frac{2|\n du|^2}{(b-\varphi\circ u)^2} + \frac{2e(u)|\n(\varphi\circ u)|^2}{(b-\varphi\circ u)^4}.
\endaligned
\end{equation}
Combining (\ref{eqn-Pa6}) with (\ref{eqn-Pa7}),  it follows
\begin{equation}\label{eqn-Pa8}\aligned
\left(\D_V - \p_t \right) f \geq & 2(1-b)(b - \varphi\circ u)f^2 
+ \frac{2\la \n f, \n(\varphi\circ u) \ra}{b-\varphi\circ u}. 
\endaligned
\end{equation}

Let $\phi := \phi  (x, t)$ be a smooth cut-off function supported in $Q_{R, T}$, satisfying the following properties (see \cite{SZ}):
\begin{itemize}
    \item[(1)] $\phi =\phi (d(x, x_0), t) \equiv \phi (r, t); \phi (x, t) = 1$ in $Q_{\frac{R}{2}, \frac{T}{4}}$ and $0\leq \phi  \leq 1.$
    \item[(2)] $\phi $ is decreasing as a radial function in the spatial variables.
    \item[(3)] $\frac{|\p_r \phi |}{\phi^a}\leq \frac{C_a}{R},  \frac{|\p_r^{2}\phi |}{\phi^a} \leq \frac{C_a}{R^2}$ for $a=\frac{1}{2}, \frac{3}{4}.$
    \item[(4)] $\frac{|\p_t \phi |}{\phi^\frac{1}{2}} \leq \frac{C}{T}$.
\end{itemize}
Let $L:= -\frac{2\n  (\varphi\circ u)}{ b - \varphi\circ u}$. Then we get
\begin{equation}\label{eqn-Pa9}\aligned
&\left(\D_V- \p_t \right)(\phi f) + \la L, \n(\phi f) \ra -2 \left\la \frac{\n \phi}{\phi}, \n(\phi f) \right\ra \\
=& f(\D_V-\p_t)\phi + \phi (\D_V-\p_t) f + 2\la \n \phi, \n f \ra \\
&+ \phi\la L, \n f \ra + f\la L, \n \phi \ra - \frac{2|\n \phi|^2 f}{\phi} - 2\la \n \phi, \n f \ra \\
=& f(\D_V-\p_t)\phi + \phi (\D_V-\p_t) f -\frac{2\phi \la \n f, \n (\varphi\circ u)\ra}{ b -\varphi \circ u } \\
& - \frac{2f \la \n \phi, \n (\varphi \circ u) \ra}{ b-\varphi\circ u } - \frac{2|\n \phi|^2f}{\phi}.
\endaligned
\end{equation}
From (\ref{eqn-Pa8}) and (\ref{eqn-Pa9}), we obtain
\begin{equation}\label{eqn-Pa10}\aligned
&\left(\D_V- \p_t \right)(\phi f) + \la L, \n(\phi f) \ra -2 \left\la \frac{\n \phi}{\phi}, \n(\phi f) \right\ra \\
\geq & 2(1-b)(b - \varphi\circ u)\phi f^2 + f \left( \D_V-\p_t \right)\phi \\
& - \frac{2f \la \n \phi, \n(\varphi\circ u) \ra}{b-\varphi\circ u} - \frac{2|\n \phi|^2f}{\phi}.
\endaligned
\end{equation}

Suppose that $\phi f$ attains its maximum at $(x_1, t_1)$.  By \cite{LY}, we can assume, without loss of generality, that $x_1$ is not in the cut-locus of $M$. Then at $(x_1, t_1)$, we have
\[
\n(\phi f) = 0, \quad \D_V(\phi f) \leq 0, \quad \p_t(\phi f) \geq 0.
\]
Then from (\ref{eqn-Pa10}), we derive at $(x_1, t_1)$,
\begin{equation}\label{eqn-para8}\aligned
2(1-b)(b - \varphi\circ u)\phi f^2 \leq  - f(\D_V-\p_t)\phi 
+ \frac{2f \la \n \phi, \n(\varphi\circ u) \ra}{b-\varphi\circ u} + \frac{2|\n \phi|^2f}{\phi}.
\endaligned
\end{equation}

Now we give an estimate for  $\D_V r$.  For any $v \in C^3 (M)$, we can derive the following Bochner type formula (see formula(2.2) in \cite{CheQiu16}):
\begin{equation}\aligned\label{V-Bochner}
\frac{1}{2}\D_V |\n v|^2 = |{\rm Hess}(v)|^2 + {\rm Ric}_V (\n v, \n v) + \la \n\D_V v, \n
v \ra.
\endaligned
\end{equation}
If $x$ is not on the cut locus of $x_0$, and for $r \geq \tilde{r}_0$($\tilde{r}_0$ is a positive constant), let $\sigma: [0, r] \rightarrow M$ be a minimal unit
speed geodesic with $\sigma(0) =x_0, \sigma(r) = x$.
Set $\varphi_V(s) = (\D_V r)\circ \sigma(s), s \in (0, r]$. Applying the Bochner formula (\ref{V-Bochner}), we have
\[
|{\rm Hess}(r)|^2  + {\rm Ric}_V(\n r, \n r) + \la \n\D_V r, \n r
\ra = \frac{1}{2}\D_V|\n r|^2 =0.
\]
Therefore,
\[
\la \n\D_V r, \n r \ra \leq -{\rm Ric}_V(\n r, \n r)
\]
Computing both sides of the above inequality along
$\sigma(s)$  gives
\begin{equation}\label{2.4}
 \varphi'_V(s) \leq - {\rm Ric}_V(\sigma', \sigma') \leq 0\
\ \ {\rm on} \ \ (0, r].
\end{equation}
Choosing $r_0>0$ small enough that geodesic ball $B_{r_0}(x_0)
\subset$ the neighborhood of the normal coordinate at $x_0$. Let
$C_1= \max_{\partial B_{r_0}(x_0)}\D_V r$. Then integrate
(\ref{2.4}) from $r_0$ to $r(x)$, we have
\begin{equation}\label{2.5}
\D_V r(x) \leq C_1.
\end{equation}
Note that
\begin{equation*}\aligned
(\D_V - \p_t )\phi = \p^2_{r} \phi |\n r|^2 +\p_r \phi \D r +\la V, \p_r \phi \n r \ra - \p_t \phi =  \p^2_{r} \phi + \p_r \phi \D_V r -\p_t \phi.
\endaligned
\end{equation*}
By (\ref{2.5}), we then obtain
\begin{equation}\label{eqn-para9}\aligned
-f(\D_V - \p_t )\phi =&  -f \p^2_{r} \phi -f  \p_r \phi \D_V r +f \p_t \phi \\
\leq & -f \p_r^{2} \phi - C_1 f\p_r \phi +f \p_t \phi \\
\leq & f |\p_r^{2}\phi| + f |C_1\p_r \phi| + f|\p_t\phi| \\
= & \phi^{\frac{1}{2}}f \cdot \frac{|\p_r^{2} \phi|}{\phi^{\frac{1}{2}}} + \phi^{\frac{1}{2}}f \cdot \frac{|C_1\p_r \phi|}{\phi^{\frac{1}{2}}} +  \phi^{\frac{1}{2}}f \cdot \frac{|\p_t \phi|}{\phi^{\frac{1}{2}}} \\
\leq & \frac{\e}{5} \phi f^2 + \frac{5}{4\e} \frac{1}{R^4} +  \frac{\e}{5}\phi f^2 \\
&+ \frac{5}{4\e} \frac{C_1^{2}}{R^2} +  \frac{\e}{5} \phi f^2 + \frac{5}{4\e} \frac{1}{T^2} \\
\leq &  \frac{3\e}{5}\e \phi f^2 + \frac{5}{4\e} \frac{1}{R^4}  + \frac{5C_1^{2}}{4\e}\frac{1}{R^2} + \frac{5}{4\e} \frac{1}{T^2}.
\endaligned
\end{equation}
Direct computation gives us
\begin{equation}\label{eqn-Pa11}\aligned
\frac{2f \la \n \phi, \n(\varphi\circ u) \ra}{b-\varphi\circ u} \leq & \frac{2f|\n \phi|\cdot e(u)^{\frac{1}{2}}}{b-\varphi\circ u} =2f^{\frac{3}{2}}|\n \phi| \\
 = &2f^{\frac{3}{2}}\phi^{\frac{3}{4}} \cdot \frac{|\n \phi|}{\phi^{\frac{3}{4}}} \leq \frac{\e}{5}\phi f^2 +\frac{3375}{16\e^3} \frac{1}{R^4}.
\endaligned
\end{equation}

The last term of the right hand side of (\ref{eqn-para8}) can be estimated as follows
\begin{equation}\label{eqn-Para11}\aligned
\frac{2|\n \phi|^2f}{\phi} = f\phi^{\frac{1}{2}}\cdot \frac{2|\n \phi|^2}{\phi^{\frac{3}{2}}} \leq \frac{\e}{5} \phi f^2 +\frac{5}{\e} \left( \frac{|\n\phi|}{\phi^{\frac{3}{4}}} \right)^4 \leq \frac{\e}{5} \phi f^2 +\frac{5}{\e} \frac{1}{R^4}.
\endaligned
\end{equation}
Substituting (\ref{eqn-para9}),  (\ref{eqn-Pa11}) and (\ref{eqn-Para11}) into (\ref{eqn-para8}), we have
\begin{equation*}\label{eqn-Para12}\aligned
2(1-b)(b - \varphi\circ u(x_1,  t_1)) \phi f^2 \leq  \e \phi f^2 +\left( \frac{25}{4\e}+ \frac{3375}{16\e^3} \right) \frac{1}{R^4} + \frac{5C_1^{2}}{4\e} \frac{1}{R^2} + \frac{5}{4\e}\frac{1}{T^2}.
\endaligned
\end{equation*}

Choosing 
\[
\e=(1-b)(b - \varphi\circ u(x_1,  t_1)).
\]
Then we get at $(x_1, t_1)$,
\begin{equation}\label{eqn-Pa12}\aligned
 \phi f^2 \leq  C \left( \left( \frac{1}{\e^2}+ \frac{1}{\e^4}  \right) \frac{1}{R^4} + \frac{1}{\e^2} \frac{1}{R^2} + \frac{1}{\e^2} \frac{1}{T^2} \right),
\endaligned
\end{equation}
where $C$ is a positive constant which is independent of $R$ and $T$.

Since 
\begin{equation}\label{eqn-QH0}
\lim_{t\to (\frac{\pi}{2})^-}\frac{\sec t}{\left( \frac{\pi}{2} -t \right)^{-1}} = \lim_{t\to (\frac{\pi}{2})^-}\frac{\frac{\pi}{2} -t}{\cos t} =  \lim_{t\to (\frac{\pi}{2})^-} \frac{-1}{-\sin t} = 1.
\end{equation}
Therefore we obtain
\[
\frac{1}{\e}=\frac{1}{(1-b)(b - \varphi\circ u(x_1,  t_1)) } \leq 4 \sup_{Q_{R, T}} \sec^2(\rho\circ u) \leq C \left( \frac{\pi}{2} - \rho\circ u \right)^{-2}.
\]
Hence from (\ref{eqn-Pa12}),  we derive
\begin{equation}\label{eqn-Pa131}
\aligned
 \sup_{Q_{\frac{R}{2}, \frac{T}{2}}} \frac{e(u)}{(b - \varphi\circ u)^{2}} 
\leq & C \left( \frac{1}{R^2}\cdot \left( \sup_{Q_{R, T}}\left( \frac{\pi}{2} - \rho\circ u \right)^{-2} 
+  \sup_{Q_{R, T}}\left( \frac{\pi}{2} - \rho\circ u \right)^{-4} \right) \right. \\
& \left.  +\frac{1}{R}\cdot  \sup_{Q_{R, T}}\left( \frac{\pi}{2} - \rho\circ u \right)^{-2} +\frac{1}{T}\cdot  \sup_{Q_{R, T}}\left( \frac{\pi}{2} - \rho\circ u \right)^{-2} \right).
\endaligned
\end{equation} 
\qed

\noindent{\bf Proof of Theorem \ref{thm2}} By (\ref{eqn-Pa131}), we have
\begin{equation}\label{eqn-Pa13}
\aligned
\frac{1}{b^{2}}\sup_{Q_{\frac{R}{2}, \frac{R^2}{2}}} e(u)
\leq & C \left( \frac{1}{R^2}\left(\sup_{Q_{R, R^2}}  \left( \frac{\pi}{2} - \rho\circ u \right)^{-2} \right.\right.  \\
& \left.\left.  + \sup_{Q_{R, R^2}}  \left( \frac{\pi}{2} - \rho\circ u \right)^{-4} \right) 
+\frac{1}{R} \sup_{Q_{R, R^2}} \left( \frac{\pi}{2} - \rho\circ u \right)^{-2} \right).
\endaligned
\end{equation}
Letting $R \to +\infty$ in (\ref{eqn-Pa13}), then we derive $e(u)\equiv 0$. Namely, $u$ is a constant.
\qed

\vskip8pt

\noindent{\bf Proof of Theorem \ref{thm-QH}}
Without loss of generality, we may assume that $\k=1$.
Let $B_R:=\{ x\in \mathbb R^{m}: |x|\leq R \}, V:=-\frac{1}{2}x$ and
\[
h_2 := (\sec\rho)^{\frac{1}{8\a}}.
\]
Direct computation gives us
\begin{equation*}\aligned
h_2'=&\frac{1}{8\a}(\sec\rho)^{\frac{1}{8\a}}\tan\rho, \\
h_2''=&(\frac{1}{8\a})^2(\sec\rho)^{\frac{1}{8\a}}\tan^2\rho + \frac{1}{8\a}(\sec\rho)^{\frac{1}{8\a}+2}.
\endaligned
\end{equation*}
Then by (\ref{eqn-Para00}), we obtain 
\begin{equation}\label{eqn-QH2}\aligned
{\rm Hess}(h_2)=& h_2'{\rm Hess}(\rho) + h_2'' d\rho\otimes d\rho \\
\geq & \frac{1}{8\a}(\sec\rho)^{\frac{1}{8\a}}\tan\rho \cot\rho (h-d\rho\otimes d\rho) \\
&+  \left( (\frac{1}{8\a})^2(\sec\rho)^{\frac{1}{8\a}}\tan^2\rho + \frac{1}{8\a}(\sec\rho)^{\frac{1}{8\a}+2} \right)d\rho\otimes d\rho \\
=&  \frac{1}{8\a}(\sec\rho)^{\frac{1}{8\a}} h +   \left( (\frac{1}{8\a})^2(\sec\rho)^{\frac{1}{8\a}}\tan^2\rho + \frac{1}{8\a}(\sec\rho)^{\frac{1}{8\a}}\tan^2\rho \right)d\rho\otimes d\rho \\
= & \frac{1}{8\a}(\sec\rho)^{\frac{1}{8\a}} h + \frac{  \frac{1}{8\a} (  \frac{1}{8\a}+1)(\sec\rho)^{\frac{1}{8\a}}\tan^2\rho }{(  \frac{1}{8\a})^2 (\sec\rho)^{\frac{1}{4\a}}\tan^2\rho} dh_2\otimes dh_2 \\
=&  \frac{1}{8\a} h_2 h + (1+8\a)h_2^{-1}dh_2\otimes dh_2.
\endaligned
\end{equation}
Denote $\widetilde{h}_2:= h_2\circ u$. By (\ref{eqn-QH2}) and the composition formula,
\begin{equation}\label{eqn-QH3}\aligned
\D_V {\widetilde h_2} = \D_V(h_2\circ u) = &\sum_{i=1}^m {\rm Hess} (h_2) (du(e_i), du(e_i))\circ u + dh_2(\tau(u)+du(V)) \\
=& \sum_{i=1}^m {\rm Hess} (h_2) (du(e_i), du(e_i)) \circ u \\
\geq & \frac{1}{8\a}{\widetilde h}_2e(u) + (1+8\a){\widetilde h}_2^{-1}|\n \widetilde h_2|^2.  
\endaligned
\end{equation}
Note that
\begin{equation}\label{eqn-QH4}
\D_V \widetilde h_2 = e^{\frac{|x|^2}{4}}\div (e^{-\frac{|x|^2}{4}}\n \widetilde h_2).
\end{equation}
From (\ref{eqn-QH3}) and (\ref{eqn-QH4}), we get
\begin{equation}\label{eqn-QH5}
\div (e^{-\frac{|x|^2}{4}}\n \widetilde h_2) \geq \frac{1}{8\a}e^{-\frac{|x|^2}{4}}\widetilde h_2 e(u).
\end{equation}
Let $\varphi_2 \in C^\infty(\mathbb R^m)$ such that $\varphi_2\equiv 1$ on $B_R, \varphi_2\equiv 0$ outside $B_{2R}$ and $|\n \varphi_2 |\leq \frac{C}{R}$. Multiplying $\varphi_2^{2}\widetilde h_2$ on both sides of (\ref{eqn-QH5}) and integrating by parts yield
\begin{equation}\label{eqn-QH6}\aligned
\frac{1}{8\a}\int_{\mathbb R^m} e(u)\varphi_{2}^2\widetilde h_2^{2} e^{-\frac{|x|^2}{4}} \leq & \int_{\mathbb R^m}\varphi_{2}^2 \widetilde h_2 \div (e^{-\frac{|x|^2}{4}} \n \widetilde h_2) = - \int_{\mathbb R^m}  \la \n(\varphi_{2}^2 \widetilde h_2), \n\widetilde h_2 \ra e^{-\frac{|x|^2}{4}} \\
=& - \int_{\mathbb R^m} 2\widetilde h_2 \varphi_{2} \la \n \varphi_2, \n \widetilde h_2 \ra e^{-\frac{|x|^2}{4}} - \int_{\mathbb R^m} \varphi_{2}^2|\n \widetilde h_2|^2 e^{-\frac{|x|^2}{4}}.
\endaligned
\end{equation}
The Cauchy inequality implies 
\begin{equation}\label{eqn-QH7}\aligned
- \int_{\mathbb R^m} 2\widetilde h_2 \varphi_{2} \la \n \varphi_{2}, \n \widetilde h_2 \ra e^{-\frac{|x|^2}{4}} \leq \frac{1}{2} \int_{\mathbb R^m} \varphi_{2}^2|\n \widetilde h_2|^2 e^{-\frac{|x|^2}{4}} + 2 \int_{\mathbb R^m} \widetilde h_2^{2} |\n \varphi_{2}|^2 e^{-\frac{|x|^2}{4}}.
\endaligned
\end{equation}
From (\ref{eqn-QH6}) and (\ref{eqn-QH7}), we obtain 
\begin{equation}\label{eqn-QH8}\aligned
\frac{1}{8\a}\int_{\mathbb R^m} e(u)\varphi_{2}^2\widetilde h_2^{2} e^{-\frac{|x|^2}{4}} +  \frac{1}{2} \int_{\mathbb R^m} \varphi_{2}^2|\n \widetilde h_2|^2 e^{-\frac{|x|^2}{4}} \leq 2 \int_{\mathbb R^m} \widetilde h_2^{2} |\n \varphi_{2}|^2 e^{-\frac{|x|^2}{4}}. 
\endaligned
\end{equation}
 It follows that
 \begin{equation*}\label{eqn-QH9}\aligned
 \int_{B_R} |\n \widetilde h_2|^2 e^{-\frac{|x|^2}{4}} \leq  \frac{4C}{R^2}\int_{B_{2R} \backslash B_R} \widetilde h_2^{2} e^{-\frac{|x|^2}{4}} \leq \frac{4C}{R^2}\sup_{B_{2R}}\widetilde h_2^{2} e^{\frac{-R^2}{4}}{\rm vol} (B_{2R}\backslash B_R).
\endaligned
\end{equation*}
 By (\ref{eqn-QH0}) and the assumption, we derive
  \begin{equation*}\label{eqn-QH10}\aligned
 \int_{B_R} |\n \widetilde h_2|^2 e^{-\frac{|x|^2}{4}} \leq & \frac{4C}{R^2}\sup_{B_{2R}}\left( \frac{\pi}{2}-\rho\circ u \right)^{-\frac{1}{4\a}}e^{\frac{-R^2}{4}}{\rm vol} (B_{2R}\backslash B_R) \\
 \leq & \frac{4C}{R^2}\sup_{B_{2R}} \left( e^{\a(\frac{1}{4}-\e_0)|x|^2} \right)^{\frac{1}{4\a}}e^{\frac{-R^2}{4}}C_m R^m \\
 \leq & C_m R^{m-2} e^{-\e_0 R^2},
\endaligned
\end{equation*}
 where $C_m$ is a positive constant depending only on $m$.
 
 Letting $R\to \infty$ forces $\widetilde h_2 \equiv c_0$, where $c_0>0$ is a constant. Again using (\ref{eqn-QH8}) gives
 \begin{equation*}\label{eqn-QH11}\aligned
\frac{c_0^{2}}{8\a}\int_{B_R} e(u) e^{-\frac{|x|^2}{4}} \leq \frac{2C}{R^2}\int_{B_{2R} \backslash B_R} \widetilde h_2^{2} e^{-\frac{|x|^2}{4}} \leq C_m R^{m-2} e^{-\e_0 R^2}.
\endaligned
\end{equation*}
Let $R\to \infty$, then we have $e(u)\equiv 0$. Hence $u$ must be a constant map.
\qed

\vskip12pt

\noindent{\bf Acknowledgements} This work is partially supported by NSFC (Nos. 11971358, 11771339) and Hubei Provincial Natural Science Foundation of China (No. 2021CFB400).

\vskip24pt

\end{document}